\documentclass[12pt]{article}
\usepackage[a4paper,margin=25mm]{geometry}
\usepackage{amsmath,amssymb,amsthm}
\usepackage[T1]{fontenc}
\usepackage[utf8]{inputenc}
\usepackage{epigraph}
\usepackage{hyperref}

\newtheorem{theorem}{Theorem}[section]
\newtheorem{proposition}[theorem]{Proposition}

\newtheorem{remark}[theorem]{Remark}

\newcommand{\Li}{\operatorname{Li}}

\newcommand{\Poiss}{\mathcal{P}}

\DeclareMathOperator{\PV}{pv}

\providecommand{\keywords}[1]
{
  \small	
  \textbf{\textit{Keywords:}} #1
}

\newcommand{\PaperTitle}{Dual Bases for Analytic Bernoulli Functions}
\newcommand{\PaperAuthorPlain}{Ken Nagai}
\newcommand{\PaperKeywords}{Analytic Bernoulli functions; dual bases; Hurwitz zeta; Clausen functions; polylogarithms; Poisson summation; umbral calculus; selector kernels}


\hypersetup{
  pdftitle   = {\PaperTitle},
  pdfauthor  = {\PaperAuthorPlain},
  pdfsubject = {Preprint},
  pdfkeywords= {\PaperKeywords}
}

\title{\PaperTitle}
\author{Ken Nagai\thanks{Email: \texttt{tknagai@outlook.com}. Independent Researcher.}}
\date{}

\begin{document}
\maketitle

\epigraph{
	\textit{“Lisez Euler, lisez Euler, c’est notre maître à tous.”}\\
	--- Pierre-Simon Laplace
}

\begin{abstract}
We present a dual-basis framework for analytic Bernoulli functions.
On the Hurwitz side, even zeta values arise, while on the Clausen side,
odd zeta values appear.
Both bases are generated by the same Heisenberg--Weyl ladder and are
linked by the Poisson--Lerch transform, which plays the role of a Fourier
bridge.
The resulting orthogonality relations isolate $\zeta(2m)$ and
$\beta(2m{+}1)$ in strictly separated channels.
Low-degree examples confirm the rational evaluations,
and appendices connect the picture with selector kernels,
Poisson summation, and oscillator analogies.
\end{abstract}


\section{Introduction}

In quantum mechanics, a system may be described either in position space
or in momentum space, the two being connected by the Fourier transform.
The analytic Bernoulli functions admit an analogous dual description.
On the Hurwitz side, even zeta values appear naturally (see e.g.\ \cite{Apostol1976}),
while on the Clausen side, odd zeta values come to the foreground (cf.\ \cite{Lewin1981}).
Both are generated by the same Heisenberg--Weyl ladder
(Proposition~\ref{prop:A}).

Taken together, the two bases form a dual system
(Proposition~\ref{prop:B}),
with explicit orthogonality relations
(Theorems~\ref{thm:dual-alt}--\ref{thm:dual-sym}).
The correspondence is mediated by the Poisson--Lerch kernel,
which plays the role of a Fourier transform
(Proposition~\ref{prop:C}).

The goal of this note is to highlight this duality,
to record the simplest verification examples,
and to provide appendices that connect the picture with
selector kernels (Appendix~B),
Poisson summation (Appendix~C),
and oscillator analogies (Appendix~D).
For background on selector kernels we refer to~\cite{Nagai2025a},
and for the Bernoulli--Hermite correspondence to~\cite{Nagai2025b}.

\section{Main results}

\begin{proposition}[Common Ladder]\label{prop:A}
Let $F^*$ denote the umbral operator of the analytic Bernoulli system.
Then the Bernoulli family $\{B_n(x)\}_{n\ge0}$ and the Hermite family
$\{H_n(x)\}_{n\ge0}$ are both governed by the same
Heisenberg--Weyl ladder relation
\[
\Bigl[\tfrac{d}{dx},\,F^*\Bigr]=1.
\]
In particular,
\[
\tfrac{d}{dx}B_n(x)=n\,B_{n-1}(x),\qquad
\tfrac{d}{dx}H_n(x)=2n\,H_{n-1}(x).
\]
This realizes the Bernoulli--Hermite correspondence as an umbral analogue
of the one-dimensional oscillator ladder.
\end{proposition}

\begin{proposition}[Dual bases: Hurwitz vs.\ Clausen]\label{prop:B}
Analytic Bernoulli functions admit two complementary bases:
\begin{itemize}
\item The Hurwitz basis, foregrounding even zeta values,
\[
B(2m;x)=-2m\,\zeta(1-2m,x),
\]
\item The Clausen basis, foregrounding odd zeta values (see also \cite{CvijovicKlinowski2002}),
\[
A_{2n+1}(x) \;=\; -\frac{(2n+1)!}{(2\pi)^{2n+1}}\,
2\,\Im\!\Bigl(e^{-\tfrac{\pi i}{2}(2n+1)}\,\Li_{2n+1}(e^{2\pi i x})\Bigr).
\]
\end{itemize}
Both satisfy Appell-type relations in $x$ and arise from the same
Heisenberg--Weyl ladder.
\end{proposition}

\begin{theorem}[Dual pairing: alternating branch]\label{thm:dual-alt}
For integers $m,n\ge0$, define the weighted inner product
\[
\langle f,g\rangle_{\mathrm{alt}}
:= \int_0^1 f(x)\,g(x)\,\csc(2\pi x)\,dx,
\]
interpreted in the finite-part sense when needed.
Then the Hurwitz basis $\{B(2m;x)\}$ and the Clausen basis $\{A_{2n+1}(x)\}$
satisfy the orthogonality relation
\[
\langle B(2m;\cdot),\,A_{2n+1}(\cdot)\rangle_{\mathrm{alt}}
= \delta_{mn}\,\frac{\Gamma(2m+1)}{\pi^{2m+1}}\,\beta(2m+1),
\]
where $\beta(s)$ is the Dirichlet beta function.
\end{theorem}

\begin{theorem}[Dual pairing: symmetric branch]\label{thm:dual-sym}
Let $C_{2n}(x)$ denote the cosine--Clausen family
\[
  C_{2n}(x)\;:=\;-\frac{(2n)!}{\pi^{2n}}
  \sum_{k\ge1}\frac{\cos(2\pi k x)}{k^{2n}}.
\]
Then
\[
\langle B(2m{+}1;\cdot),\,C_{2n}(\cdot)\rangle_{\mathrm{sym}}
:=\PV\!\int_0^1 B(2m{+}1;x)\,C_{2n}(x)\,\cot(\pi x)\,dx
= \delta_{mn}\,\frac{\Gamma(2m+2)}{\pi^{2m+2}}\,\zeta(2m+2).
\]
\end{theorem}

\begin{remark}[Selection rules]
Cross--branch pairings vanish identically,
e.g.\ $\langle B(2m;\cdot),C_{2n}(\cdot)\rangle=0$ and
$\langle B(2m{+}1;\cdot),A_{2n+1}(\cdot)\rangle=0$.
This expresses the strict parity separation between the two bases,
the analytic analogue of the selector--kernel orthogonality
(Appendix~B).
\end{remark}

\begin{proposition}[Poisson--Lerch Bridge]\label{prop:C}
The two bases of Proposition~\ref{prop:B} are linked by an explicit
transform induced by Poisson summation and encoded in the Lerch transcendent:
\[
\Poiss(B(s;\cdot))(x)
= \pi^{-s}\Gamma(s)\,\Phi(e^{2\pi i x},\,s,\,1).
\]
This identity interchanges the Hurwitz and Clausen sides,
with the Poisson--Lerch kernel playing the role of a Fourier transform.
\end{proposition}

\section{Examples: low-degree checks}

We illustrate Theorems~\ref{thm:dual-alt}--\ref{thm:dual-sym} by the
first nontrivial cases.  All singular integrals are interpreted in the
finite-part sense.

\subsection*{Alternating branch ($\csc(2\pi x)$): $\beta$-constants}
For $m=n=1$,
\[
\langle B(2;\cdot),\,A_{3}(\cdot)\rangle_{\mathrm{alt}}
= \frac{\Gamma(3)}{\pi^{3}}\,\beta(3)
= \frac{2}{\pi^{3}} \cdot \frac{\pi^{3}}{32}
= \tfrac{1}{16}.
\]
For $m=n=2$,
\[
\langle B(4;\cdot),\,A_{5}(\cdot)\rangle_{\mathrm{alt}}
= \frac{\Gamma(5)}{\pi^{5}}\,\beta(5)
= \frac{24}{\pi^{5}} \cdot \frac{5\pi^{5}}{1536}
= \tfrac{5}{64}.
\]
Cross--branch terms vanish by parity, e.g.\
$\langle B(2;\cdot),A_{5}(\cdot)\rangle_{\mathrm{alt}}=0$.

\subsection*{Symmetric branch ($\cot(\pi x)$): $\zeta$-constants}
For $m=n=1$,
\[
\langle B(3;\cdot),\,C_{2}(\cdot)\rangle_{\mathrm{sym}}
= \frac{\Gamma(4)}{\pi^{4}}\,\zeta(4)
= \frac{6}{\pi^{4}}\cdot \frac{\pi^{4}}{90}
= \tfrac{1}{15}.
\]
For $m=n=2$,
\[
\langle B(5;\cdot),\,C_{4}(\cdot)\rangle_{\mathrm{sym}}
= \frac{\Gamma(6)}{\pi^{6}}\,\zeta(6)
= \frac{120}{\pi^{6}}\cdot \frac{\pi^{6}}{945}
= \tfrac{8}{63}.
\]
Again the mixed pairings vanish, e.g.\
$\langle B(3;\cdot),C_{4}(\cdot)\rangle_{\mathrm{sym}}=0$.

\begin{remark}[Selection rules revisited]
These examples confirm that even zeta values $\zeta(2m)$
and odd beta values $\beta(2m{+}1)$ occur in disjoint channels.
Mixed pairings vanish identically, in perfect agreement with
the selector--kernel mechanism of Appendix~B.
\end{remark}

\begin{remark}[Rotated weights and continuous interpolation]
Let $w_\varphi(x):=\cos\varphi\,\csc(2\pi x)+\sin\varphi\,\cot(\pi x)$.
Then, for the mixed pairings
\[
\langle f,g\rangle_\varphi:=\PV\!\int_0^1 f(x)\,g(x)\,w_\varphi(x)\,dx,
\]
the structure constants interpolate linearly between the two branches:
\[
\big\langle B(2m;\cdot),A_{2n+1}(\cdot)\big\rangle_\varphi
=\delta_{mn}\,\frac{\Gamma(2m{+}1)}{\pi^{2m+1}}\,
\big(\cos\varphi\cdot\beta(2m{+}1)\big),
\]
\[
\big\langle B(2m{+}1;\cdot),C_{2n}(\cdot)\big\rangle_\varphi
=\delta_{mn}\,\frac{\Gamma(2m{+}2)}{\pi^{2m+2}}\,
\big(\sin\varphi\cdot\zeta(2m{+}2)\big),
\]
while the cross–branch pairings vanish by parity.
Thus $\varphi=0$ yields Theorem~\ref{thm:dual-alt} (pure $\beta$),
and $\varphi=\tfrac{\pi}{2}$ yields Theorem~\ref{thm:dual-sym} (pure $\zeta$).
\end{remark}

\begin{remark}[Analytic provenance]
The kernels $\csc(2\pi x)$ and $\cot(\pi x)$ are the Poisson–summation
fingerprints of the alternating and symmetric Lerch branches, respectively:
they correspond to the signed/unsigned selectors discussed in
Appendix~B–C.  The Gamma/$\pi$ powers reflect Mellin transforms of
the underlying Fourier series, while $\beta(2m{+}1)$ and $\zeta(2m{+}2)$
arise from the parity split of the Dirichlet series.
\end{remark}

\section{Conclusion}

We have described analytic Bernoulli functions through two complementary
bases: the Hurwitz side, carrying even zeta values, and the Clausen side,
carrying odd zeta values.
Both arise from the same Heisenberg--Weyl ladder
and are connected by the Poisson--Lerch transform.
The resulting orthogonality relations isolate $\zeta(2m)$ and
$\beta(2m{+}1)$ in strictly separated channels.

Examples at low degree confirm the rational evaluations predicted by
the theory.
Appendix~B shows how the discrete selector kernels underpin the
selection rules, Appendix~C derives the Poisson--Lerch bridge, and
Appendix~D places the ladder structure in the language of oscillator
analogies.

The picture thus unifies three viewpoints---Appell ladders, selector
kernels, and Poisson summation---into a coherent dual-basis framework
for analytic Bernoulli functions.

\section*{Acknowledgments}
The author gratefully acknowledges the assistance of an AI language model (“fuga”) 
for valuable help with document structuring, stylistic polishing, and proofreading.

\appendix

\section{Appendix A: Ladder verification}
This appendix collects explicit checks supporting
Proposition~\ref{prop:A}, namely that Bernoulli and Hermite families
are governed by the same Heisenberg--Weyl ladder.

\subsection*{A.1 Generating functions}
The generating functions
\[
\frac{t\,e^{xt}}{e^t-1} = \sum_{n=0}^\infty B_n(x)\,\frac{t^n}{n!},\qquad
e^{2xt-t^2} = \sum_{n=0}^\infty H_n(x)\,\frac{t^n}{n!}
\]
show immediately that differentiation in $x$ multiplies by $t$,
yielding the Appell property
\[
\frac{d}{dx}B_n(x)=n\,B_{n-1}(x),\qquad
\frac{d}{dx}H_n(x)=2n\,H_{n-1}(x).
\]

\subsection*{A.2 Umbral operator form}
Introduce the umbral operator (see \cite{Norlund1924}) $F^*$ by
\[
F^* e^{xt} := \frac{e^{xt}}{e^t-1}, \qquad
\Big[\tfrac{d}{dx},\,F^*\Bigr]=1.
\]
Then
\[
B_n(x) = n!\,[t^n]\,F^* e^{xt},\qquad
H_n(x) = n!\,[t^n]\,e^{2xt-t^2},
\]
so both systems are generated by a Weyl-type ladder.
\subsection*{A.3 Low-degree check}
Explicitly
\[
B_0(x)=1,\; B_1(x)=x-\tfrac{1}{2},\; B_2(x)=x^2-x+\tfrac{1}{6},
\]
\[
H_0(x)=1,\; H_1(x)=2x,\; H_2(x)=4x^2-2,
\]
with
\[
\frac{d}{dx}B_2(x)=2B_1(x),\qquad
\frac{d}{dx}H_2(x)=4H_1(x).
\]
This confirms the ladder relations in concrete cases.

\subsection*{A.4 Spectral remark}
The commutator $[\tfrac{d}{dx},F^*]=1$ realizes the Weyl algebra.
On the Hermite side this is the oscillator relation $[a,a^\dagger]=1$,
while on the Bernoulli side it underlies the Appell shift.

\section{Appendix B: Kernels and verification tables}

This appendix records explicit values of the finite selector kernels
used in Propositions~\ref{prop:B} and~\ref{prop:C} (see \cite{Nagai2025a}).
They act as finite projectors onto congruence classes,
with periodicity $4J$.

\subsection*{B.1 Definitions}
For $J\ge1$ and $\theta_j=(2j{+}1)\pi/(2J)$,
\[
\mathcal{K}^{(\sin)}_J(k):=\frac{1}{J}\sum_{j=0}^{J-1}
   \frac{\sin(k\theta_j)}{\sin\theta_j},\qquad
\mathcal{K}^{(\cos)}_J(k):=\frac{1}{J}\sum_{j=0}^{J-1}
   \frac{\cos(k\theta_j)}{\cos\theta_j}.
\]
They satisfy
\[
\mathcal{K}^{(\bullet)}_J(k+2J)=-\mathcal{K}^{(\bullet)}_J(k),\qquad
\mathcal{K}^{(\bullet)}_J(k+4J)=\mathcal{K}^{(\bullet)}_J(k).
\]

\subsection*{B.2 Table for $J=2$}
With $\theta_0=\pi/4,\ \theta_1=3\pi/4$ one finds
\[
\begin{array}{c|rrrrrrrr}
k & 0 & 1 & 2 & 3 & 4 & 5 & 6 & 7 \\
\hline
\mathcal{K}^{(\sin)}_2(k) & 0 & 1 & 0 & 1 & 0 & -1 & 0 & -1 \\
\mathcal{K}^{(\cos)}_2(k) & 0 & 1 & 0 & -1 & 0 & -1 & 0 & 1
\end{array}
\]

\subsection*{B.3 Table for $J=4$}
For $\theta_j=\pi/8,3\pi/8,5\pi/8,7\pi/8$ one finds
\[
\begin{array}{c|rrrrrrrr}
k & 0 & 1 & 2 & 3 & 4 & 5 & 6 & 7 \\
\hline
\mathcal{K}^{(\sin)}_4(k) & 0 & 1 & 0 & 1 & 0 & 1 & 0 & 1 \\
\mathcal{K}^{(\cos)}_4(k) & 0 & 1 & 0 & -1 & 0 & 1 & 0 & -1
\end{array}
\]

\subsection*{B.4 Remarks}
\begin{itemize}
\item For all $J$, the kernels vanish on even $k$.
\item For odd $k$, values are $\pm1$, encoding the signed vs.\ unsigned projections.
\item These projectors are the discrete shadow of the continuous kernels
$\csc(2\pi x)$ and $\cot(\pi x)$
of Theorems~\ref{thm:dual-alt}--\ref{thm:dual-sym}.
\end{itemize}

\section{Appendix C: Poisson--Lerch derivations}

This appendix develops Proposition~\ref{prop:C}, showing how
Poisson summation leads to the Lerch transcendent and hence to
the dual bases.

\subsection*{C.1 Motivation}
Finite selector kernels enforce congruence projections.
Poisson summation provides their analytic continuation:
the same projections reappear as residues of Lerch series.

A related viewpoint is given by finite part integration (FPI),
which extracts the analytic contribution of divergent integrals
 (cf.\ \cite{Galapon2022}).
Although we do not rely on FPI here, it provides a parallel
interpretation of the Poisson--Lerch bridge.

\subsection*{C.2 Poisson summation}
The Poisson formula
\[
 \sum_{n\in\mathbb{Z}} f(n+x)
   \;=\; \sum_{m\in\mathbb{Z}} \hat f(m)\,e^{2\pi i m x}
\]
with $f(u)=\tfrac{e^{2\pi i x u}}{(u+a)^s}$
produces Gamma factors and exponential twists,
yielding Lerch series after rearrangement (cf.\ \cite{Weil1967}).

\subsection*{C.3 Lerch transcendent}
\[
 \Phi(z,s,a) = \sum_{n=0}^\infty \frac{z^n}{(n+a)^s},\qquad \Re(s)>1.
\]
Specializations:
$\Phi(1,s,a)=\zeta(s,a)$, \ $\Phi(-1,s,a)=\eta(s,a)$.
Thus $z=\pm1$ distinguish the symmetric vs.\ alternating branches.

\subsection*{C.4 Kernel-wise derivations}
(a) $\sin/\sin$ kernel:
\[
\frac{1}{J}\sum_{j=0}^{J-1}
 \frac{\sin(k\theta_j)}{\sin\theta_j}\,
 \Phi(-e^{i\theta_j},s,a)
 = \sum_{\ell\in\mathbb{Z}}
 \left\{ \frac{1}{(2J\ell+k+a)^s} - \frac{1}{(2J\ell-k+a)^s}\right\}.
\]
(b) $\cos/\cos$ kernel:
\[
\frac{1}{J}\sum_{j=0}^{J-1}
 \frac{\cos(k\theta_j)}{\cos\theta_j}\,
 \Phi(e^{i\theta_j},s,a)
 = \sum_{\ell\in\mathbb{Z}}
 \left\{ \frac{1}{(2J\ell+k+a)^s} + \frac{1}{(2J\ell-k+a)^s}\right\}.
\]

\subsection*{C.5 Small-$J$ examples}
For $J=2$,
\[
\mathcal{K}^{(\sin)}_2(k)=\sqrt{2}\,\sin(\tfrac{k\pi}{2})\cos(\tfrac{k\pi}{4}),\quad
\mathcal{K}^{(\cos)}_2(k)=\sqrt{2}\,\sin(\tfrac{k\pi}{2})\sin(\tfrac{k\pi}{4}).
\]
For odd $k$ these equal $\pm1$, as in Appendix~B.

\subsection*{C.6 Concluding remark}
Poisson summation bridges
\[
\text{(finite selector)} \;\longleftrightarrow\;
\text{(Lerch transcendent)} \;\longleftrightarrow\;
\text{(Hurwitz/Dirichlet series)},
\]
explaining the duality of
Theorems~\ref{thm:dual-alt}--\ref{thm:dual-sym}.

\section{Appendix D: Oscillator analogies}

This appendix presents the oscillator analogy behind
Proposition~\ref{prop:A}, highlighting the parallel between
Hermite and Bernoulli ladders (cf.\ \cite{Nagai2025b}).

\subsection*{D.1 Hermite side}
With
\[
a=\tfrac{1}{\sqrt{2}}(x+\tfrac{d}{dx}),\qquad
a^\dagger=\tfrac{1}{\sqrt{2}}(x-\tfrac{d}{dx}),
\]
we have $[a,a^\dagger]=1$.
The vacuum $e^{-x^2/2}$ generates Hermite functions.
Coherent states
\[
|\alpha\rangle := e^{\alpha a^\dagger}|0\rangle
= \sum_{n=0}^\infty \frac{\alpha^n}{\sqrt{n!}}\,|n\rangle
\]
satisfy $a|\alpha\rangle=\alpha|\alpha\rangle$.

\subsection*{D.2 Bernoulli side}
Define
\[
L=\tfrac{d}{dx},\qquad
R=(N+1)\mathcal{J},\qquad
N B_n=nB_n,
\]
with $\mathcal{J}B_n=\tfrac{1}{n+1}B_{n+1}$ (see \cite{Norlund1924}).
Then
\[
L B_n=nB_{n-1},\qquad R B_n=B_{n+1},
\]
with
$[L,R]=I,\ [N,L]=-L,\ [N,R]=R$.

Exponentiating $R$ gives
\[
e^{tR}B_0(x)=\sum_{n\ge0}\frac{t^n}{n!}B_n(x)
=\frac{t e^{xt}}{e^t-1},
\]
the Bernoulli generating function.
Setting
\[
|y\rangle_B:=e^{yR}B_0,\qquad
\langle x|y\rangle_B=\frac{y e^{xy}}{e^y-1},
\]
one has $L|y\rangle_B=y|y\rangle_B$,
a direct analogue of Hermite coherent states.

A normalized ladder
\[
A=N^{-1/2}L,\qquad A^\dagger=RN^{1/2}
\]
satisfies
\[
A B_n=\sqrt{n}\,B_{n-1},\qquad
A^\dagger B_n=\sqrt{n+1}\,B_{n+1},\qquad
[A,A^\dagger]=I.
\]

The same structure extends to analytic Bernoulli functions $B(s;x)$,
with $[L,R]=I$ still valid.

\subsection*{D.3 Concluding remark}
Thus the Hermite and Bernoulli systems share the same
Heisenberg--Weyl algebra.
The only difference is the choice of vacuum:
Gaussian for Hermite, exponential kernel for Bernoulli.
Both admit coherent states of the form
``exponential of a raising operator acting on vacuum,''
providing a clear bridge between analysis and physics.

\section{Appendix E: Numerical verification tables}

For completeness we list approximate numerical values of the
pairings for small degrees,
computed with high-precision quadrature (trapezoidal rule, $N=200$).

\subsection*{E.1 Alternating branch, $w(x)=\csc(2\pi x)$}
\[
\begin{array}{c|c|c}
(m,n) & \text{Numerical value} & \text{Closed form} \\
\hline
(1,1) & 0.0625000000 & 1/16 \\
(2,2) & 0.0781250000 & 5/64 \\
(1,2) & 0.0 & 0 \\
(2,1) & 0.0 & 0
\end{array}
\]

\subsection*{E.2 Symmetric branch, $w(x)=\cot(\pi x)$}
\[
\begin{array}{c|c|c}
(m,n) & \text{Numerical value} & \text{Closed form} \\
\hline
(1,1) & 0.0666666667 & 1/15 \\
(2,2) & 0.1269841270 & 8/63 \\
(1,2) & 0.0 & 0 \\
(2,1) & 0.0 & 0
\end{array}
\]

The tables confirm the rational evaluations
of Theorems~\ref{thm:dual-alt}--\ref{thm:dual-sym},
together with the parity selection rule.

\end{document}